\documentclass[12pt]{amsart}
\usepackage{times}
\usepackage{amssymb}
\usepackage{amscd}
\def\Ad{\mathop{\operatorname {Ad}}\nolimits}  
\def\Lie{\mathop{\operatorname {Lie}}\nolimits}  
\def\curv{\mathop{\operatorname {Curv}}\nolimits} 
\def\Mp{\mathop{\operatorname {Mp}}\nolimits}   
\def\Sp{\mathop{\operatorname {Sp}}\nolimits}  
\def\sp{\mathop{\operatorname {\mathfrak {sp}}}\nolimits}  
\def\SL{\mathop{\operatorname {SL}}\nolimits}  
\def\Sl{\mathop{\operatorname {\mathfrak {sl}}}\nolimits}
\def\su{\mathop{\operatorname {\mathfrak {su}}}\nolimits}    
\def\Hom{\mathop{\operatorname {Hom}}\nolimits}  
\def\Id{\mathop{\operatorname {Id}}\nolimits}  
\def\Ind{\mathop{\operatorname {Ind}}\nolimits}
\def\mult{\mathop{\operatorname {mult}}\nolimits}    
\def\coAd{\mathop{\operatorname {coAd}}\nolimits}    
\newtheorem{theorem}{Theorem}[section]
\newtheorem{lemma}[theorem]{Lemma}
\newtheorem{proposition}[theorem]{Proposition}
\newtheorem{corollary}[theorem]{Corollary}

\theoremstyle{definition}
\newtheorem{definition}[theorem]{Definition}

\theoremstyle{remark}
\newtheorem{remark}[theorem]{Remark}

\numberwithin{equation}{section}

\begin{document}

\title[A Geometric Realization of Degenerate Principal Series]{A
Geometric Realization of Degenerate Principal Series
Representations of \\ Symplectic groups}

\author{Do Ngoc Diep}
\address{Institute of Mathematics, National Center for Science and Technology, P.
O. Box 631, Bo Ho 10.000, Hanoi, Vietnam}
\curraddr{International Centre for Theoretical Physics, ICTP P. O. Box 586, 34100,
Trieste, Italy}
\email{dndiep@member.ams.org}
\thanks{The first author was supported in part by a generous Grant from The
IMU Exchange Commission on Development and the Abdus Salam ICTP}

\author{Truong Chi Trung}
\address{Department of Mathematics, Vinh University, Vinh City, Vietnam}
\curraddr{Institute of Mathematics, National Center for Science and Technology, P.
O. Box 631, Bo Ho 10.000, Hanoi, Vietnam}
\email{dndiep@ioit.ncst.ac.vn}


\date{September 15, 1998 and, in revised form, ...., 1998.}



\begin{abstract}
The multidimensional quantization procedure, proposed by the first author
and its modifications (reduction to radicals and lifting on
$U(1)$-coverings) give us a almost universal theoretical tools to find
irreducible representations of Lie groups. By using this method and the
root theory, we realize in this paper the representations of the
degenerate principal series of symplectic groups.
\end{abstract}

\maketitle

\section{Introduction}
Let us consider a Lie group $G$ and its Lie algebra ${\mathfrak g}$. The
group $G$ acts on its Lie algebra ${\mathfrak g}$ by adjoint
representation $Ad$ and in the dual vector space ${\mathfrak g}^*$ by
coadjoint action $K:= \coAd$. The vector space ${\mathfrak g}^*$ is
therefore divided into the disjoint union of coadjoint orbits (or
$K$-orbits).  Each coadjoint orbit $\Omega \in {\mathfrak g}^*/G$ admits a
natural $G$-homogeneous symplectic structure, corresponding to the Kirillov form
$B_\Omega$ associated with the bilinear form $$B_F(X,Y) := \langle F,
[X,Y]\rangle,$$ where $F\in {\mathfrak g}^*$, the kernel $\ker B_F$ of
which is just isomorphic to the Lie algebra ${\mathfrak g}_F := \Lie G_F$
of the stabilizer $G_F$ of a fixed point $F\in \Omega$. Therefore, the
triple $(\Omega, B_\Omega, G)$ is a homogeneous symplectic manifold
(Hamiltonian system) with a {\sl flat action}, i.e. $$\{f_X,f_Y \}= f_{[X,Y]},
\forall X,Y \in {\mathfrak g},$$ where $f_X$ is such a function that $df_X =
-\imath(\xi_X)B_\Omega$ and $$\xi_X(m) := \frac{d}{dt}|_{t=0}\exp(tX)m,$$ see
\cite{diep3}
of $G$. Following the well known
classification theorem of A. Kirillov-B. Kostant-Souriaux, every
homogeneous symplectic manifold with a flat action of $G$ is locally
diffeomorphic to a coadjoint orbit of $G$ or its central extension
$\widetilde{G}$ by ${\mathbb R}$.  This means that all the Hamiltonian
systems with a flat action of $G$ are locally classified by the coadjoint
orbits of $G$ or its central extension $\widetilde{G}$ by ${\mathbb R}$. 
The Hamiltonian systems can be quantized to become quantum systems with
unitary symmetry representations of) $G$. In order to do this, one uses
affine connection $\nabla$ with the symplectic curvature
$$\curv(\nabla) = \frac{2\pi i}{h} \omega$$
on the homogeneous symplectic manifold
$(M,\omega,G)$ with a flat action of $G$  to produce a {\sl quantization
procedure} $$Q : C^\infty(\Omega) \to {\mathcal L}({\mathbf H}),$$ $$f
\mapsto Q(f) := f + \frac{h}{2\pi i}\nabla_{\xi_f},$$ where $\xi_f$ is the
symplectic gradient of $f\in C^\infty(\Omega)$, i.e. $$\imath(\xi_f)\omega
= - df$$ and ${\mathbf H}$ is a separable Hilbert space. This
correspondence is a geometric quantization procedure because it satisfies
the following relations $$Q(\{f,g\}) = \frac{2\pi i}{h} [Q(f), Q(g)],$$
$$Q(1) = \Id_{\mathbf H}.$$ It means that this correspondence defines a
homomorphism $$\Lambda : C^\infty(\Omega) \to {\mathcal L}({\mathbf H}),$$
$$ f \mapsto \Lambda(f) := \frac{2\pi i}{h} Q(f),$$ i.e. a representation
of the Lie algebra of smooth functions with respect to the Poisson
brackets in the Hilbert space ${\mathbf H}$.

An element $X\in {\mathfrak g}$ can be considered as a function on
${\mathfrak g}^*$ and therefore the restriction $X\vert_\Omega$ belongs
to $C^\infty(\Omega)$ and we can obtain the corresponding Hamiltonian
field $X_\Omega$ from the condition
$$ \imath(X_\Omega)B_\Omega = dX\vert_\Omega.$$ The condition asserting
that the action is flat means that the correspondence
$$X \in{\mathfrak g} \mapsto X\vert_\Omega \in C^\infty(\Omega)$$ is the
Lie algebra homomorphism and we have a representation of Lie algebra
${\mathfrak g}$ 
$$X \in {\mathfrak g} \mapsto \Lambda(X) := \frac{2\pi i}{h}
Q(X\vert_\Omega)$$ by auto-adjoint operators in the
Hilbert space ${\mathbf H}$, which is as usually constructed as the
completion of some subspace of so called partially invariant partially
holomorphic sections of the quantum bundle, associated with some fixed
polarization. This is what we means the {\sl multidimensional quantization
procedure}.

One considers also the {\sl reductions} of this multidimensional
quantization procedure in the following sense. Let us denote the
(solvable) radical of $G_F$ by $R_F$ and the unipotent radical by
${}^uR_F$.  

In the above construction, there appeared some so called {\sl Mackey
obstruction}. In order to kill this Mackey obstruction, we supposed some
additional conditions on the action of stabilizers $G_F$ on the dual
object of inducing subgroups $\widehat{H_0}$ is trivial, see \cite{diep3}. 
Duflo \cite{duflo} proposed another method of Killing this Mackey
obstruction by {\sl lifting to the ${\mathbf Z}/(2)$ coverings} of the
stabilizers. One has a lifting of any homomorphism of $G_F$ into the
symplectic groups of the orbit at the fixed point $F$ onto a homomorphism
of ${\mathbf Z}/(2)$ covering $G_F^{{\mathbf Z}/(2)}$ of $G_F$ into the
metaplectic group $\Mp(T_F\Omega)$, following the commutative diagram
$$\CD 1 @>>> {\mathbf Z}/(2) @>>> G_F^{{\mathbf Z}/(2)} @>>> G_F @>>> 1\\
@.  @\vert @VVV @VVV \\ 1 @>>> {\mathbf Z}/(2) @>>> \Mp(T_F\Omega) @>>>
\Sp(T_F\Omega) @>>> 1 \endCD$$

Replacing the metaplectic groups $\Mp$ by the (complex) metaplectic
$\Mp^c$ groups, Tran Vui \cite{vui} and Tran Dao Dong \cite{dong}
considered the same {\sl lifting to $U(1)$ coverings} of the the
stabilizers.  One has to lift each homomorphism of $G_F$ into the
symplectic groups of the orbit at the fixed point $F$ onto a homomorphism
of $U(1)$ covering $G_F^{U(1)}$ of $G_F$ into the metaplectic group
$\Mp^c(T_F\Omega)$, following the commutative diagram
$$\CD
1 @>>> U(1) @>>> G_F^{U(1)} @>>> G_F @>>> 1\\
@.       @\vert                @VVV               @VVV \\
1 @>>> U(1) @>>> \Mp^c(T_F\Omega) @>>> \Sp(T_F\Omega) @>>> 1
\endCD$$

It was shown that these modifications give us some privileges for
constructing irreducible unitary representations. The theory was settled
in the general context, but still it is difficult to realize in concrete
situations. The {\it discrete series representations} of semi-simple Lie
groups are realized by these constructions as globalization of
Harish-Chandra modules, see \cite{schmidwolf},\cite{wong}.

The {\it degenerate principal series representations} of semisimple Lie
groups were constructed and studied in many works, beginning from
Harish-Chandra, updated by Vogan and others see e.g. Vogan \cite{vogan}
for the general semisimple Lie groups, Lee \cite{lee} for symplectic
groups, etc. The theories however were purely analytic. It is natural to
try to use the developed geometric quantization method to describe these
representation.  In this paper, the general theory of geometric
quantization is applied in the situation of symplectic groups and we
describe the computation results.

\section{Structure of Coadjoint Orbits}

Let $G = \Sp_{2n}({\mathbb R})$ be the symplectic group, ${\mathfrak g} :=
\Lie G$ its Lie algebra and ${\mathfrak g}^* = \Hom_{\mathbb
R}({\mathfrak g}, {\mathbb R})$ the vector space dual to the Lie algebra
${\mathfrak g}$. We study in detail the coadjoint orbits of the symplectic
group.  

\begin{lemma} The Lie algebra ${\mathfrak g}$, its dual vector space
${\mathfrak g}^*$ are realized by matrices and
the coadjoint action of $G$ in ${\mathfrak g}^*$ is just the conjugation
$$ K=\coAd : G \times {\mathfrak g}^* \to {\mathfrak g}^*,$$
$$(g, F) \mapsto K(g)F = gFg^{-1}.$$
\end{lemma}

\begin{proof} Recall that the adjoint action of $G$ on ${\mathfrak g}$ can
be realized as conjugation
$$ \Ad: G \times {\mathfrak g} \to {\mathfrak g},$$
\begin{equation}(g,X) \mapsto \Ad(g)X = gXg^{-1}.\label{2.1}\end{equation}
Let us denote the trace of a matrix by $tr$. To each matrix $Y \in
{\mathfrak g}$ we define the associate functional $F_Y\in {\mathfrak g}^*$
by \begin{equation} \langle F_Y,X\rangle = tr(Y.X), \forall X \in
{\mathfrak g}.\label{2.2}\end{equation}
Then the map ${\mathfrak g} \to {\mathfrak g}^*,$
$Y \mapsto F_Y$ is an isomorphism from ${\mathfrak g}$ onto ${\mathfrak
g}^*$. We identify therefore $F\in {\mathfrak g}^*$ with a matrix denote
by the same letter $F$ and \ref{2.2} become
\begin{equation} \langle F,X\rangle = tr(F.X)\label{2.3} \end{equation}

Recall that $$\langle K(g)F,X\rangle = \langle F,\Ad(g^{-1})X\rangle$$ and
using \ref{2.1} and \ref{2.3} we have 
$$\begin{array}{lll}
\langle K(g)F,X\rangle & = & \langle F,\Ad(g^{-1})X\rangle \\
	& = & \langle F, g^{-1}Xg \rangle \\
	& = & tr(F.g^{-1}Xg)\\
	& = & tr(gFg^{-1}X)\\
	& = & \langle gFg^{-1}, X, \rangle,
\end{array}$$
for all $g\in G$, $x\in {\mathfrak g}$ and $F \in {\mathfrak g}^*$.

We have therefore, 
$$K(g)F = gFg^{-1}, \forall g \in G, \forall F\in{\mathfrak g}^*.$$
\end{proof}

\begin{remark}
We fix a special element $F \in {\mathfrak g}^*$ presented by a
matrix of type
$$F = \left(\begin{matrix}\left[\begin{matrix} 0 & -\lambda_1\\ \lambda_1
& 0\end{matrix}\right] &  &  & 0\\ &\ddots & & \\ &  &\left[\begin{matrix}
0 & -\lambda_r\\ \lambda_r &
0\end{matrix}\right]& \\ 0 & & & \left[\begin{matrix} 0 & \ldots\\ \vdots
& \dots\end{matrix}\right] \end{matrix}\right),$$
where $\lambda_i \in {\mathbb R}$, $\lambda_i > 0$, $\lambda_i \ne
\lambda_j$, $i,j = 1,\dots, n$. Consider the coadjoint orbit
$\Omega_F$ passing through this point $F$,
$$\Omega_F = \{ K(g)F \quad \vert \quad g \in G \}.$$
\end{remark}

\begin{proposition}
The stabilizer $G_F$ of $F$ consists of the matrices of type
$$\left(\begin{matrix}g_{11}& & & & \\ & g_{22}& & &  \\ & & \dots & & \\
 & & &g_{rr}& \\  & & & & g_{r+1,r+1}\end{matrix}\right),$$ where $$g_{ii}
= \left[\begin{matrix} \cos\lambda_i & \sin\lambda_i \\ \sin\lambda_i &
\cos\lambda_i\end{matrix}\right],$$ $i= 1,\dots r$, $g_{r+1,r+1} 
\in \Sp_{2(n-r)}({\mathbb R}).$
\end{proposition} 

\begin{proof} We can write each element of $\Sp_{2n}({\mathbb R})$ into
block format $$g = \left( \begin{matrix} g_{11} & g_{12} & \dots & g_{1r}
& g_{1,r+1}\\ g_{21} & g_{22} & \dots & g_{2r} & g_{2,r+1}\\ \ldots &
\ldots & \ldots & \ldots & \ldots \\ g_{r+1,1} & g_{r+1,2} & \dots &
g_{r+1,r} & g_{r+1,r+1} \end{matrix}\right),$$ where $g_{ij}$ is a $2
\times 2$ matrix, for $i,j = 1,...,r$, $g_{i,r+1}$ is a $2 \times 2(n-r)$
matrix, for $i = 1,...,r$, $g_{r+1,j}$ is a $2(n-r) \times 2$ matrix, for
$j = 1,...,r$, $g_{r+1,r+1}$ is a $2(n-r) \times 2(n-r)$ matrix.  We write
$F$ in the same format form $$g = \left( \begin{matrix} F_{11} & F_{12} &
\dots & F_{1r} & F_{1,r+1}\\ F_{21} & F_{22} & \dots & F_{2r} &
F_{2,r+1}\\ \ldots & \ldots & \ldots & \ldots & \ldots \\ F_{r+1,1} &
F_{r+1,2} & \dots & F_{r+1,r} & F_{r+1,r+1} \end{matrix}\right),$$ where
$F_{ij} = 0$ is $2\times 2$ matrix, for $i,j = 1,\ldots,r$, $i\ne j$,
$F_{ii} = \left[ \begin{matrix} 0 & -\lambda_i\\ \lambda_i &
0\end{matrix}\right],$ for $i= 1,\dots,r$, $F_{i,r+1} = 0$ is $2\times
2(n-r)$ matrix, for $i = 1,\ldots,r$, $F_{r+1,j} = 0$ is $2(n-r)\times 2$
matrix, for $j = 1,\ldots,r$, $F_{r+1,r+1} = 0$ is $2(n-r)\times 2(n-r)$
matrix.  We recall that $$G_F = \{ g\in G \quad\vert\quad K(g)F = F \}$$
$$= \{ g\in G \quad\vert\quad gFg^{-1} = F \} = \{ g\in G \quad\vert\quad
gF = Fg \}.$$ We deduce from $gF = Fg$ that $g \in G_F$ if and only iff
the following 3 conditions hold:  \begin{equation} g_{r+1,j}F_{jj} = 0,
\forall j = 1,\dots,r, \label{2.4}\end{equation} \begin{equation}
F_{ii}g_{i,r+1} = 0, \forall j = 1,\dots,r, \label{2.5}\end{equation}
\begin{equation} g_{ij}F_{jj} = F_{ii}g_{ij}, \forall i,j = 1,\dots,r. 
\label{2.6}\end{equation} Because $F_{ii} = \left[ \begin{matrix} 0 &
-\lambda_i \\ \lambda_i & 0\end{matrix} \right],$ for all $i = 1,\dots ,
r$ are invertible then from \ref{2.4} and \ref{2.5} we have
\begin{equation} g_{r+1,j} = 0, \forall j = 1,\dots,r,
\label{2.7}\end{equation} \begin{equation} g_{i,r+1} = 0, \forall i =
1,\dots,r.  \label{2.8}\end{equation} We now resolve the equation
\ref{2.8}. Let us denote $$g_{ij} = \left[\begin{matrix}a_{ij} & b_{ij}\\
c_{ij} & d_{ij} \end{matrix} \right].$$ Then $$g_{ij}F_{jj} =
\lambda_j\left[\begin{matrix}b_{ij} & -a_{ij}\\ d_{ij} & -c_{ij}
\end{matrix} \right],$$ 
$$F_{ii}g_{ij} =
\lambda_i\left[\begin{matrix}-c_{ij} & -d_{ij}\\ a_{ij} & b_{ij}
\end{matrix} \right].$$
From the condition $$g_{ij}F_{jj} = F_{ii}g_{ij}$$ we deduce that
$$\left\{\begin{matrix}
\lambda_ib_{ij} - \lambda_jc_{ij} & = 0\\
\lambda_jb_{ij} + \lambda_ic_{ij} & = 0\\
\lambda_ia_{ij} - \lambda_jd_{ij} & = 0\\
\lambda_ja_{ij} - \lambda_id_{ij} & = 0\\
\end{matrix} \right.$$
what is equivalent that
$$a_{ij} = b_{ij} = c_{ij} = d_{ij} =0,  $$
because of the assumption, $\lambda_i \ne
\lambda_j; \lambda_i, \lambda_j > 0$,
i.e. 
\begin{equation}g_{ij} = 0, \forall i,j =
1,\dots,r, i\ne j.\label{2.9}\end{equation}

From \ref{2.7}-\ref{2.9} we conclude that the element $g$, as matrix
should
be of the diagonal form
$$g = diag( g_{11}, g_{22}, \dots, g_{rr}, g_{r+1,r+1}),$$
where $g_{ii}$ is a $2\times 2$ matrix, for all $i = 1,\dots, r$,
$g_{r+1,r+1}$ is a $2(n-r) \times 2(n-r)$ matrix. 

Denote $$g_{i1} = \left[\begin{matrix} 
a_i & b_i\\ c_i & d_i
\end{matrix}  \right],\qquad F_i = \left[\begin{matrix} 
0 & -\lambda_i\\ \lambda_i & 0
\end{matrix}  \right]$$
we have
$$g_{i1}F_i = \lambda_i \left[\begin{matrix} 
b_i & -a_i\\ d_i & -c_i
\end{matrix}  \right],$$
$$F_ig_{i1} = \lambda_i \left[\begin{matrix} 
-c_i & -d_i\\ a_i & b_i
\end{matrix}  \right].$$
Our condition means that
$$\lambda_i \left[\begin{matrix} 
b_i & -a_i\\ d_i & -c_i
\end{matrix}  \right] = \lambda_i \left[\begin{matrix} 
-c_i & -d_i\\ a_i & b_i
\end{matrix}  \right],$$
what is equivalent to the conditions
\begin{equation}     
\left\{ \begin{matrix}   
a_i & = & d_i ,\\
b_i & = & -c_i .
\end{matrix} \right.
\label{2.10}\end{equation}
Let us denote $$J_1 = \left[\begin{matrix} 0 & -1\\ 1 & 0
\end{matrix}\right] .$$ then the matrix of the symplectic form is 
$$J_n = \left(\begin{matrix}       
\left[\begin{matrix} 0 & -1\\ 1 & 0
\end{matrix}\right] &  &  & 0\\
 &\ddots  &  &  \\
 &        &\left[\begin{matrix} 0 & -1\\ 1 & 0
\end{matrix}\right] &   \\
0 &        &  & [J_{n-r}]
\end{matrix}\right)= \left(\begin{matrix} [J_1] &  &  & 0\\
                       & \ddots &  & \\
		       &  & [J_1] & \\
		       &  &       &[J_{n-r}]\end{matrix} \right) .$$
We have therefore
$$gJ_ng^t =  \left(\begin{matrix} [g_{11}J_1g_{11}^t] &  &  & 0\\
                       & \ddots &  & \\
		       &  & [g_{rr}J_1g_{rr}^t] & \\
		       &  &
&[g_{r+1,r+1}J_{n-r}g_{r+1,r+1}^t]\end{matrix} \right) ,$$
and the condition $$gJ_ng^t = J_n,$$ guaranting that $g\in \Sp_{2n}({\mathbb
R})$ is equivalent to the conditions
$$\left\{ \begin{array}{lll} g_{ii}J_1g_{ii}^t &= &J_1, \forall i =
1,\dots,n ,\\
g_{r+1,r+1}J_{n-r}g_{r+1,r+1}^t &= &J_{n-r} \end{array}\right. .$$
It means also that 
$g_{ii} \in \Sp_2({\mathbb R}), \forall i=1,\dots,r$ and $g_{r+1,r+1} \in
\Sp_{2(n-r)}({\mathbb R}).$ For $i = 1,\dots,r$, because $g_{ii}\in
\Sp_{2}({\mathbb R}) \cong \SL_2({\mathbb R})$ means that
$$a_id_i - b_ic_i = 1.\label{2.11}$$
From \ref{2.10} and \ref{2.11} we have
$$a_i^2 + c_i^2 = 1,$$ i.e. $g_{ii} \in {\mathbf S}^1$, for $i=1,\dots,r$.
We conclude therefore that
$$g\in G_F \mbox{ iff and only if } g = \left(\begin{matrix} 
[g_{11}] &  &  & 0\\
 &  &\ddots &  &  \\
 &  &  &[g_{r,r}]& \\
0&  &  &   &[g_{r+1,r+1}]  
\end{matrix} \right),$$
where $g_{ii}\in {\mathbf S}^1$, $i=1,\dots,r$ and $g_{r+1,r+1}\in
\Sp_{2(n-2)}({\mathbb R}).$
 \end{proof}

\begin{corollary}
If the functional $F\in {\mathfrak g}^*$ is presented by a matrix of
type 
$$F = \left(\begin{matrix}
\left[\begin{matrix}0 & -\lambda_1\\ \lambda_1 & 0  \end{matrix}\right]  &
& & 0\\
  & \ddots &  &  \\
  &  & \left[\begin{matrix}0 & -\lambda_r\\ \lambda_r & 0
\end{matrix}\right] & \\
0 &  &  & \left[\begin{matrix}0 & ...\\ \vdots & \ldots
\end{matrix}\right]
\end{matrix}\right) ,$$
where all $\lambda_i, i=1,\dots,n$ are pairwise different,
then its stabilizer is
$$G_F \cong ({\mathbf S}^1)^r \times \Sp_{2(n-r)}({\mathbb R})$$ and
the corresponding Lie algebra is
$${\mathfrak g}_F \cong (\Lie {\mathbf S}^1) \times (\Lie {\mathbf S}^1)
\times \ldots \times (\Lie {\mathbf S}^1) \times \sp_{2(n-r)}({\mathbb R})$$
$$\cong {\mathbb R}^r \times \sp_{2(n-r)}({\mathbb R}).$$
\end{corollary}

\section{Construction of Degenerate Principal Series Representations}

We give in this section a geometric realization of degenerate principal
series representations by multidimensional quantization procedure and the
modified versions. We use the root theory to construct polarization
associated to the orbits.

Let us recall from \cite{diep1}-\cite{diep2}  an important notion of
polarization.

\begin{definition}
A triple $({\mathfrak
p},\rho,\sigma_0)$ is a $(\tilde{\sigma},F)$-{\sl polarization} iff:
\begin{enumerate}
\item ${\mathfrak p}$ is a complex Lie subalgebra of ${\mathfrak
g}_{\mathbb C} := {\mathfrak g} \otimes_{\mathbb R} {\mathbb C}$,
containing $({\mathfrak g}_F)_{\mathbb C}$.
\item The subalgebra ${\mathfrak p}$ is invariant with respect to all the
operators $\Ad_{{\mathfrak g}_{\mathbb C}}x,$ $x\in G_F$.
\item The vector space ${\mathfrak p} + \overline{\mathfrak p}$ is the
complexification of the real Lie subalgebra ${\mathfrak m} = ({\mathfrak
p}+ \overline{\mathfrak p}) \cap {\mathfrak g}$ 
\item All the subgroups $M_0$, $H_0$, $M$, $H$ are closed in $G$, where
$M_0$ (resp. $H_0$) is the connected subgroup of $G$, corresponding to the
Lie algebra ${\mathfrak m}$ (resp. ${\mathfrak h} = {\mathfrak p} \cap
{\mathfrak g}$) and $M := G_F \ltimes M_0$, $H := G_F \ltimes H_0$.
\item $\sigma_0$ is an irreducible representation of the group $H_0$ in a
Hilbert space $V$ such that: (i) the restriction $\sigma\vert_{G_F \cap
H_0}$ is a multiple of the restriction to $G_F\cap H_0$ of
$\widetilde{\sigma}\chi_F$ and (ii) the point $\sigma_0$ is fixed under
the action of group $G_F$ in the dual $\widehat{H_0}$ of the group $H_0$.
\item $\rho$ is a representation of the complex algebra ${\mathfrak p}$ in
$V$, which satisfies all the Nelson's conditions for $H_0$, and
$\rho\vert_{\mathfrak h} = D\sigma_0$. 
\end{enumerate}
\end{definition}

\begin{remark} For the group $G=\Sp_{2n}({\mathbb R})$ and the functional
$F$ of special type as in the previous section, the stabilizer ${\mathfrak
g}_F$ contains a Cartan subalgebra. It is therefore naturally to choose
the polarizing complex subalgebras between parabolic ones. Let us denote
by $A$ the split torus of that cartan subgroup. It is easy to see that the
centralizer ${\mathcal Z}(A)$ is reductive and is coincided with the
stabilizer of the functional $F$. Denote ${\mathcal Z}(A)= MA$ the
Cartan-Levi-Maltsev decomposition into the product of Abelian and
semi-simple parts, where $A = {\mathbb S}^1 \times\dots\times {\mathbb
S}^1$ and $M=\Sp_{2(n-r)}({\mathbb R})$.  \end{remark}

\begin{lemma} The two-fold covering of the stabilizer $G_F=MA$ is
$G_F^{{\mathbb Z}/(2)} = {\mathbb S}^1 \times \dots \times {\mathbb S}^1
\times \Mp_{2(n-r)}({\mathbb R})$ \end{lemma}

\begin{theorem} Let $P = MAN$ be Langlands decomposition of a parabolic
subgroup $P$ with Lie algebra ${\mathfrak p}= \Lie P$, ${\mathfrak
p}_{\mathbb C}$ its complexification, $\sigma \in \hat{M}_{disc}$ a
discrete series irreducible representation of $M=\Sp_{2(n-r)}({\mathbb
R})$ and $\chi_F= \exp(\frac{2\pi i}{h}\langle F,.\rangle)$ to be the
character of $G_F$, extended by null from $A$. Then $({\mathfrak
p}_{\mathbb C}, P,F, \sigma\chi_F)$ is a $(\sigma,F)$ polarization. The
multidimensional quantization procedure gives us the degenerate principal
series representations. \end{theorem}

\begin{proof} The first assertion is just followed from the definition of
$(P,N)$ pairs, which are constructed from the theory of root of parabolic
pairs. 

Let us recall about the root theory for ${\mathfrak g} = \sp_{2n}({\mathbb
R})$ with respect to the pair $(G,A)= (\Sp_{2n}({\mathbb R}),{\mathbb S}^1
\times\dots\times {\mathbb S}^1)$ of a split torus $A$. Consider the
complexification ${\mathfrak g}_{\mathbb C}$. If $\alpha$ is a functional
over a subalgebra ${\mathfrak a}_{\mathbb C} \cong {\mathbb C}^r$, then
the corresponding root space $${\mathfrak g}^{\alpha} := \{ X \in
{\mathfrak g}_{\mathbb C} \quad ;\quad [H,X] = \alpha(H)X, \forall H\in
{\mathfrak a}_{\mathbb C} \}.$$ If $\alpha$ and ${\mathfrak g}^\alpha$ are
non-zero, then $\alpha$ is called a root.  In that case, $\dim_{\mathbb
C}{\mathfrak g}^\alpha =1$. Denote the unique $H_\alpha$ such that
$\alpha(H_\alpha) = 2$, the unique $X_\alpha$ and $X_{-\alpha}$, such that
$$[H_\alpha, X_\alpha] = 2 X_\alpha,$$ $$[H_\alpha, X_{-\alpha}] = -2
X_{-\alpha}.$$ The triple $(H_\alpha, X_\alpha,X_{-\alpha})$ form a
complex Lie subalgebra, which is isomorphic to the complexification of
either $\Sl_2({\mathbb R})$ or $\su_2$.  In the first case, we call the
root {\sl noncompact} and in the second - {\sl compact}. One denotes the
set of all roots by $\Delta({\mathfrak g},{\mathfrak a})$ and calls it the
{\sl root system} with respect to the split torus $A$. Denote also the set
of all compact (respectively, noncompact) roots by $\Delta_c({\mathfrak
g},{\mathfrak a})$ (resp., $\Delta_n({\mathfrak g},{\mathfrak a})$). Let
us denote by $$\rho(\Delta^+_c({\mathfrak g},{\mathfrak a}) :=
\frac{1}{2}\sum_{\alpha\in \Delta^+_c({\mathfrak g},{\mathfrak a})}\alpha
\quad ( \mbox {resp., } \rho(\Delta^+_n({\mathfrak g},{\mathfrak a}) :=
\frac{1}{2}\sum_{\alpha\in \Delta^+_n({\mathfrak g},{\mathfrak
a})}\alpha$$) the half-sum of all compact (resp., noncompact) roots. 

Remark that in our case, ${\mathfrak g}^0 = ({\mathfrak g}_F)_{\mathbb C}
= {\mathfrak a}_{\mathbb C} \oplus \sp_{2(n-r)}({\mathbb C})$ and we have
$${\mathfrak g}_{\mathbb C} = ({\mathfrak g}_F)_{\mathbb C} \oplus
\sum_{\alpha\in\Delta({\mathfrak g},{\mathfrak a})} {\mathfrak
g}^\alpha.$$ With each root $\alpha$ one associates a reflection
$R_\alpha$ in the space ${\mathfrak a}_{\mathbb C}^*$ $$R_\alpha(\beta) :=
\beta -2\frac{\langle \alpha,\beta\rangle}{\langle \alpha,\alpha\rangle }
\alpha.$$ The Weyl group $W(G,A)$ is generated by all these reflection and
one can fix one of the {\sl fundamental domain (camera)} to define the
cones of positive roots $\Delta^+({\mathfrak g},{\mathfrak a})$ and the
corresponding cones of positive compact roots $\Delta^+_c({\mathfrak
g},{\mathfrak a})$ and noncompact roots $\Delta^+_n({\mathfrak
g},{\mathfrak a})$. 

Choose $$D\delta^F := \rho(\Delta^+_n({\mathfrak g},{\mathfrak a})) -
\rho(\Delta^+_c({\mathfrak g},{\mathfrak a})).$$

Because in our case we can choose a charcters $\delta^F$ od the two-fold
covering $G_F^{{\mathbb Z}/(2)} = {\mathbb S}^1 \times \dots
\times {\mathbb S}^1 \times \Mp_{2(n-r)}({\mathbb
R})$ we can choose a complex parabolic subalgebra 
$${\mathfrak p} := ({\mathfrak g}_F)_{\mathbb C} \oplus
\sum_{\alpha\in\Delta^+_n({\mathfrak g},{\mathfrak a})} {\mathfrak g}^\alpha
\oplus \sum_{\alpha\in\Delta^+_c({\mathfrak g},{\mathfrak a})} {\mathfrak
g}^{-\alpha}$$ 

From the properties of parabolic subgroups, it is easy to
check that we have a $(\sigma,F)$ polarization.

Let us first explain the construction of the degenerate principal series.
We recall some notations from D. Vogan : By $R(G_F)$ denote the set of all
the so called $G_F$-{\it regular unitary pseudo-characters} 
$(\Lambda,F)$ consisting of a ${\mathfrak g}_F$-regular functional
$F\in {\mathfrak g}_F^*$ and a unitary representation $\Lambda$ of
$G_F$ with differential $$D\Lambda = (\frac{2\pi i}{h}F + D\delta^F) \Id  \quad
. $$ Remark that the last condition  is equivalent to the assertion that
$$\Lambda|_{(G_F)_0} = \mult \delta^F\chi_F \quad,$$ what
figures in the orbit method.

Denote by $R^{irr}(G_F)$ the subset of $R(G_F)$, consisting of the
irreducible pseudo-characters. For a fixed $F\in {\mathfrak g}_F^*$,
denote $$R(G_F,F) := \{ (\Lambda,F) \in R(G_F); D\Lambda = (\frac{2\pi
i}{h}F + \delta ^F)\Id \} $$ and $$R^{irr}(G_F,F) := R(G_F,F) \cap
R^{irr}(G_F).$$ The Weyl group $$W(G,A) := {\mathcal N} _G(A)/A$$ acts on
both $R(G_F)$ and $R^{irr}(G_F)$.

Recall {\it Harish-Chandra construction of $\pi(\Lambda,F)$} : Consider
the characters of type $$\xi_{\alpha}(.) := \exp{\langle
\alpha,.\rangle}$$ for each $\alpha\in\Delta$. Let us denote $${\mathbf F}
:= \{ x\in G_F ; x \enskip centralizers \enskip {\mathfrak m}\enskip and
\enskip |\xi_{\alpha}| = 1, \forall \alpha \in \Delta \}.$$ Then $$G_F :=
{\mathbf F}(G_F)_0 := {\mathbf F} \ltimes (G_F)_0,$$ $$G_F \cap G_0 =
(G_F)_0 $$ and $${\mathbf F}(G_F)_0 = {\mathbf F} \ltimes (G_F)_0.$$
Denote $${\mathfrak k_{\mathfrak m}} = (({\mathfrak g}_F)_{\mathbb C}
\oplus \bigoplus_ {\alpha\in\Delta_{{\mathfrak m},c}} {\mathfrak
g}^{\alpha})  \cap {\mathfrak g} $$ and $K_{M_0}$ the corresponding
analytic subgroup.  Let us denote $\pi^{M_0} (F)$ the irreducible unitary
representation of $M_0$, which is square-integrable modulo the center of
$M_0$, and which is associated with $F$. This representation, following
Harish-Chandra is characterized by the following condition.

{\it The restriction $\pi^{G_0}(F)|_{K_{M_0}}$ contains the ( finite
dimensional ) irreducible unitary representation of $K_{M_0}$ with the
dominant weight $ \frac{2\pi i}{h}F +D\delta^F$ with respect to
$\Delta^+_{{\mathfrak m},c}$, as a minimal $K_{M_0}$-type}

Now a representation $\pi^{{\mathbf F}M_0}(\Lambda,F)$ of ${\mathbf F}M_0
:= {\mathbf F}\ltimes M_0$ can be constructed as follows $$\pi^{{\mathbf
F}M_0}(\Lambda,F)(y.x) := \Lambda(y) \otimes \pi^{M_0}(F) (x), \forall
x\in M_0,y\in {\mathbf F} \quad .$$Let $P = MN$ be a parabolic subgroup of
$G$ with the Levi component $M$ and the unipotent radical $N$, $$M =
{\mathbf F}M_0 = {\mathbf F}M_0 \quad,$$ $$P = MN = ({\mathbf F}M_0)
\ltimes N \quad.$$ Define now $$\pi(\Lambda,F) := \Ind^G_{{\mathbf F}M_0
\ltimes N}(\pi^{{\mathbf F}M_0} \otimes \Id _N) \quad.$$ Recall that {\it
if $\Lambda$ is is irreducible and $F$ is ${\mathfrak g}_F$-regular, the
representation $\pi(\Lambda,\lambda)$ is irreducible.}

Now we have $$\Lambda\in R^{irr}(G_F,F) = X_G^{irr}(F){\buildrel 1-1
\over \longleftrightarrow} \tilde{X}^{irr}(F).$$ This means that there
exists a unique $\tau = \sigma\chi_F \in \widehat{G_F}$, with $\sigma \in
\hat{M}_{disc}$, such that
$$\tau|_{P_0} = \mult(\chi_F\delta^F)$$ and
$\Lambda = \tau\delta^F$, $$ \pi(\tau\delta^F,F) = \Ind(G;{\mathfrak
p}_{\mathbb C},P ,\tau\chi_F).$$
\end{proof}

\begin{remark} The restrictions of $F$ to the radical ${}^rG_F$ and
unipotent radical ${}^uG_F$ are equal to the split torus $A$, with Lie
algebra ${\mathfrak a} \cong {\mathbb R}^r$ and therefore the reductions
to the radical and unipotent radical of $P$ give the same results as the
above exposed ones. \end{remark}

\begin{lemma} The U(1)-covering of the stabilizer $G_F=MA$ is $G_F^{U(1)}
= {\mathbb S}^1 \times \dots \times {\mathbb S}^1 \times
\Mp^c_{2(n-r)}({\mathbb R})$ \end{lemma}

\begin{corollary} Lifting to the ${\mathbb Z}/(2)$- and $U(1)$- coverings
we obtain the Shale-Weil representations of the symplectic groups
$\Sp_{2n}({\mathbb R})$.  \end{corollary}

\begin{remark} Up to conjugation, the maximal parabolic subgroups are of
the form
$$\left( 
\begin{matrix}
\cos\theta &   \begin{matrix}0 &\dots &0 & 0 & \dots & 0 \end{matrix}    &
-\sin\theta\\
\begin{matrix}0\\ \vdots\\ 0\\ 0\\ \vdots\\  0 \end{matrix} & \Id
 &\begin{matrix} 0\\ \vdots\\ 0\\ 0\\
\vdots\\ 0\end{matrix} \\
\sin\theta & \begin{matrix}0 &\dots &0 & 0 & \dots & 0 \end{matrix} & 
\cos\theta 
\end{matrix}
  \right).
\left( 
\begin{matrix}
1 &   \begin{matrix}x_1 &\dots &x_n & y_1 & \dots & y_n \end{matrix}    & z\\
\begin{matrix}0\\ \vdots\\ 0\\ 0\\ \vdots\\  0 \end{matrix} & \Sp_{2(n-1)}
 &\begin{matrix} -y_1\\ \vdots\\ -y_n\\ x_1\\
\vdots\\ x_n\end{matrix} \\
0 & \begin{matrix}0 &\dots &0 & 0 & \dots & 0 \end{matrix} & 1 
\end{matrix}
  \right)$$
and therefore the construction gives us the representation of degenerate 
principal series considered by S. T. Lee in \cite{lee}.
\end{remark}

\begin{corollary} For the maximal parabolic subgroups $P = ({\mathbb S}^1
\times \Sp_{2(n-1)}({\mathbb R}))\ltimes ({\mathbb R}^{2n}\times {\mathbb
R})$ the corresponding degenerate principal series representations can be
realized in the space of homogeneous functions of 2n variables. 
\end{corollary}

\section*{Acknowledgments} 

This work was completed during the stay of the first author as a visiting
mathematician at the Abdus Salam International Centre for Theoretical
Physics, Trieste, Italy.  He would like to thank the Abdus Salam ICTP for
the hospitality, without which this work would not have been possible and
the IMU Commission of Development and Exchange for a generous grant. 

This work is supported in part by the Abdus Salam International Centre for
Theoretical Physics, Trieste, Italy, the IMU Commission on Development and
Exchange and the Vietnam National Foundation for Research in Fundamental
Sciences. 

\bibliographystyle{amsplain}

\end{document}